\documentclass[12pt,a4paper,twoside]{article}

\usepackage{amsmath,amsfonts,amsthm,amsopn,color,amssymb}
\usepackage{graphicx}
\newcommand{\eps}{\varepsilon}

\newcommand{\R}{\mathbb{R}}

\newcommand{\RN}{{\mathbb{R}^N}}

\DeclareMathOperator{\cat}{cat}

\DeclareMathOperator{\dist}{dist}

\renewcommand{\le}{\leqslant}
\renewcommand{\ge}{\geqslant}
\renewcommand{\a }{\alpha }

\renewcommand{\b }{\beta }

\renewcommand{\d }{\delta }

\newcommand{\g }{\gamma }

\renewcommand{\l }{\lambda}
\newcommand{\n }{\nabla }

\renewcommand{\t}{\theta}

\newcommand{\G}{\Gamma}

\newcommand{\Ne}{\mathcal{N}}

\newcommand{\M}{\mathcal{M}}

\renewcommand{\o}{\omega}

\newcommand{\D }{{\mathcal D}^{1,2}(\RN)}
\newcommand{\irn }{\int_{\RN}}
\newcommand{\intrn }{\int_{\RN}}

\def\bbm[#1]{\mbox{\boldmath $#1$}}

\newtheorem{theorem}{Theorem}[section]
\newtheorem{lemma}[theorem]{Lemma}

\newtheorem{remark}[theorem]{Remark}

\renewenvironment{proof}{\noindent{\textbf{Proof\quad}}}{$\hfill\square$\vspace{0.2 cm}\\}
\newenvironment{proofmain}{\noindent{\textbf{Proof of Theorem \ref{th:main}\quad}}}{$\hfill\square$\vspace{0.2 cm}\\}

\title{{\bf On a ``zero mass'' nonlinear Schr\"odinger~equation}}

\author{A. Azzollini \thanks{Dipartimento di Matematica, Universit\`a degli
Studi di Bari,  Via E. Orabona 4, I-70125 Bari, Italy, e-mail: {\tt azzollini@dm.uniba.it}}
 \; \& \;
A. Pomponio\thanks{Dipartimento di Matematica, Politecnico di Bari, Via Amendola 126/B, I-70126 Bari, Italy, e-mail: {\tt a.pomponio@poliba.it}}}

\date{}

\begin{document}

\maketitle

\begin{abstract}
We look for positive solutions to the nonlinear Schr\"odinger~equation $-\eps^2\Delta u - V(x)f'(u)=0$, in $\RN$, where $V$ is a continuous bounded positive potential and $f$ satisfies particular growth conditions which make our problem fall in the so called ``zero mass case". We prove an existence result for any $\eps>0,$ and a multiplicity result for $\eps$ sufficiently small.
\end{abstract}

\section{Introduction and statement of the results}

In this paper we study the elliptic equation
\begin{equation}\label{eq}\tag{$\cal{P}$}
\left\{
\begin{array}{l}
-\Delta u - V(x)f'(u)=0, \qquad \hbox{ in }\RN,
\\
u>0,
\\
u \in \D,
\end{array}
\right.
\end{equation}
where $N\ge 3$, $V:\RN \to \R$ and $f:\R\to \R$. We are interested in the so called ``zero mass case" that is, roughly speaking, when $f''(0)=0$.

When $V$ is a positive constant, such a problem has been intensely studied by many authors. Some results have been obtained by \cite{Au,CGS,T}, if $f$ corresponds to the critical power $t^{(N+2)/(N-2)}$, and by \cite{BL1,BL2,BL3,P}, if $f$ is supercritical near the origin and subcritical at infinity (see also \cite{BM} for the case of exterior domain and \cite{AP} for complex valued solutions).

Up to our knowledge, there is no result in the literature on problem \eqref{eq} when $V$ is not a constant. Our aim is to investigate this case. More precisely, we will assume the following hypotheses on $V:\RN\to \R$
\begin{itemize}
\item[({\bf V1})] $V\in C(\RN,\R)$;
\item[({\bf V2})] $C_1\le V(x)\le C_2$, for all $x\in \RN$;
\item[({\bf V3})] $\limsup_{|y|\to\infty}V(y) \le V(x)$, for all $x\in \RN$, and the
inequality is strict for some $x\in \RN$;
\end{itemize}
and $f:\R\to\R$ satisfying
\begin{itemize}
\item[({\bf f1})] $f\in C^2(\R,\R)$ and even;
\item[({\bf f2})] $\forall t\in\R:$ $f(t)\ge C_3  \min(|t|^p,|t|^q)$;
\item[({\bf f3})] $\forall t\in\R:$ $|f'(t)|\le C_4 \min(|t|^{p-1},|t|^{q-1})$;
\item[({\bf f4})] $\exists \a>2$ such that $\forall t \in\R\setminus\{0\}:\;\a f(t)\le f'(t)t< f''(t)t^2$;
\end{itemize}
with $2<p<2^*=(2N)/(N-2)<q$ and $C_1,$ $C_2,$ $C_3,$ $C_4$  positive constants. These particular growth conditions on $f$ were introduced by \cite{BF} to study the semilinear Maxwell equations.

We get the following result:
\begin{theorem}\label{th:main}
If $f$ satisfies $(${\bf f1-4}$)$ and $V$ satisfies $(${\bf V1-3}$)$, then problem \eqref{eq} possesses at least a nontrivial solution.
\end{theorem}

We also consider the singularly perturbed version of problem
\eqref{eq}, namely we look for solutions of the problem
\begin{equation}\label{eq:eps}\tag{$\cal{P}_\eps$}
\left\{
\begin{array}{l}
-\eps^2 \Delta u - V(x)f'(u)=0,\qquad \hbox{ in }\RN,
\\
u>0,
\\
u \in \D,
\end{array}
\right.
\end{equation}
for $\eps>0$ sufficiently small.

Replacing ({\bf V3}) by
\begin{itemize}
\item[({\bf V4})] $\limsup_{|x|\to\infty}V(x) <\sup_{x\in \RN} V(x)$,
\end{itemize}
we get the following result:
\begin{theorem}\label{th:maineps}
If $f$ satisfies $(${\bf f1-4}$)$ and $V$ satisfies $(${\bf V1-2}$)$ and $(${\bf V4}$)$, then problem \eqref{eq:eps} possesses at least a nontrivial solution, for $\eps$ sufficiently small.
\end{theorem}
Observe that, since
({\bf V3}) implies ({\bf V4}), the introduction of a small parameter $\eps>0$
allows us to obtain an existence result assuming weaker hypotheses on the potential $V$.

Actually the introduction of the parameter $\eps$ allows us to get a stronger result then Theorem \ref{th:maineps}.

We set
\begin{align*}
M:=&\,\Big\{\eta\in \RN \;\big|\; V(\eta)=\max_{\xi\in \RN}V(\xi)\Big\},
\\
\noalign{ \hbox{and for any $\g>0$,}}
M_\g:=&\,\Big\{\eta\in\R^N\;\big|\;\inf_{\xi\in M}\|\eta-\xi\|_{\R^N}\le
\g\Big\}.
\end{align*}
Observe that $M\neq \emptyset$, by ({\bf V4}).

We get the following multiplicity result
\begin{theorem}\label{th:main2}
If $V$ satisfies $(${\bf V1-2}$)$, $(${\bf V4}$)$ and $f$
satisfies $(${\bf f1-4}$),$ then, for every $\g>0$, there exists
$\bar \eps>0$ such that the problem \eqref{eq:eps} has at least
$\cat_{M_\g}(M)$ nontrivial solutions for any $\eps \in (0,\bar
\eps)$.
\end{theorem}
Here $\cat_{M_\g}(M)$ means the
Lusternik-Schnirelmann category of $M$ in $M_\g$.

\begin{remark}
The evenness of the nonlinearity $f$ is required just to obtain positive solutions for problems \eqref{eq} and \eqref{eq:eps}. If we are not interested in the sign of solutions, the evenness hypothesis can be removed.
\end{remark}

This paper has been motivated by some well known works, such as
\cite{ABC,AMS,CL1,CL2,dPF,FW,Oh1,Oh2,R,W,WZ}, where the nonlinear Schr\"odinger equation
\[
-\eps^2 \Delta u + K(x)u = R(x)|u|^{r-2} u, \qquad \hbox{in }\RN
\]
has been studied for $2<r<2^*$ in the ``positive mass case'',
namely when $K$ is bounded below by a positive constant (see also
\cite{AF} for the p-Laplacian).

In Section \ref{sec:ex}, we take a variational approach to
\eqref{eq} and \eqref{eq:eps}. The key point of this section is
Theorem \ref{th:c}, which, in the same spirit of \cite{BN,R},
guarantees the compactness at mountain pass level under a suitable
value. By means of this theorem, we prove Theorems \ref{th:main}
and \ref{th:maineps}. Even if Theorem \ref{th:maineps} follows
immediately from Theorem \ref{th:main2}, we prefer to prove it
directly in this section since it is strictly correlated with
Theorem \ref{th:c}.

In Section \ref{sec:mult}, following \cite{AF,CL1,CL2}, we look at the topological and
compactness properties of the sublevels of the functional
associated to \eqref{eq:eps}, in order to prove Theorem
\ref{th:main2}.

\section{Existence results}\label{sec:ex}

Throughout all this section, we will suppose that the hypotheses ({\bf f1-4}) and ({\bf V1-2}) hold.

In order to find weak solutions of the
problem \eqref{eq}, we define the functional
$I\colon \D \to \R$ as:
\[
I(u)=\frac 12 \irn |\n u|^2 - \irn V(x) f(u)\,d x,
\]
where $\D$ is the completion of $C^\infty_0(\RN)$ with respect to
the norm
\[
\|u\|=\left(\irn |\n u|^2\right)^\frac 12.
\]
Observe that, by the growth condition ({\bf f3}), the functional $I$ is well defined and of class $C^1$, and its critical points correspond  to weak solutions of \eqref{eq}.
Moreover we denote by $\Ne$ the so called Nehari manifold of $I$, namely
\[
\Ne :=\left\{u \in \D \setminus \{0\} \;\Big{|}\; \irn |\n u|^2 = \irn
V(x) f'(u)u\,d x \right\}.
\]
Using similar arguments as those in \cite{BM}, we can prove
\begin{lemma}\label{le:N}
\begin{enumerate}
\item $\Ne$ is a $C^1$ manifold;
\item for any $u\neq 0$ there exists a unique number $\bar\t>0$ such that $\bar\t u\in \Ne$ and
\[
I(\bar\t u)=\max_{\t \ge 0}I(\t u);
\]
\item there exists a positive constant $C$, such that for all $u \in \Ne$, $\|u\|\ge C$.
\end{enumerate}
\end{lemma}
By 2 of Lemma \ref{le:N}, the map $\t:\D\setminus\{0\}\to\R_+$
such that for any $u\in\D,$ $u\neq 0:$
\begin{equation*}
I(\t (u) u)=\max_{\t\ge 0} I(\t u)
\end{equation*}
is well defined.
\\
Set
\begin{align}
c_1 &= \inf_{g\in \G} \max_{\t\in [0,1]} I(g(\t));
\label{eq:c1}\\
c_2 &= \inf_{u\neq 0} \max_{\t \ge 0} I(\t u);\nonumber
\\
c_3 &= \inf_{u\in\Ne} I(u);\nonumber
\end{align}
where
\[
\G=\left\{g\in C\big([0,1],\D\big) \mid g(0)=0,\;I(g(1))\le 0, \;g(1)\neq 0\right\}.
\]
Arguing as in \cite[Proposition 3.11]{R}, we also have
\begin{lemma}\label{le:ccc}
The following equalities hold
\[
c=c(V):=c_1=c_2=c_3.
\]
\end{lemma}

Observe that, since we are in unbounded domain, there is a lack of compactness. In
particular, it is in general not true that the (PS)-sequences,
namely sequences of the type $(u_n)_n\subset\D$ such that
\begin{eqnarray*}
\big(I(u_n)\big)_n \hbox{ is bounded and }
I'(u_n)\to 0,
\end{eqnarray*}
admit a converging subsequence. Moreover, the presence of the
potential $V$ does not permit us to use any symmetry to recover
compactness in a suitable natural constraint of $\D.$ In order to
overcome this difficulty, we are going to use a
concentration-compactness argument as in \cite{ABDF} (see also \cite{L1,L2}).
\\
The following lemma provides the boundedness and the concentration
of the (PS)-sequences (actually we consider a more general
situation).
\\
In the sequel $(V_n)_n$ is a sequence of potentials satisfying
({\bf V1-2}) uniformly, $(I_n)_n$ is the sequence of the
functionals defined by
    \begin{equation*}
        I_n(u):=\frac 1 2 \irn|\n u|^2-\irn V_n(x) f(u)
    \end{equation*}
and $c_n:=c(V_n).$

\begin{lemma}\label{le:tea}
Let $0<a\le b.$ If $(u_n)_n \subset \D$ is such that
    \begin{equation*}
        a\le I_n(u_n)\le b\quad\hbox{and}\quad I_n'(u_n)\to
        0,
    \end{equation*}
then
\begin{enumerate}
\item $(u_n)_n$ is bounded in the $\D$;
\item there exist a sequence $(y_n)_n \subset \RN$ and two positive numbers $R,\;\mu
>0$ such that
\begin{equation}\label{eq:br}
\liminf_n \int_{B_R(y_n)}|u_n|^2 \, d x >\mu.
\end{equation}
In particular, there exists a positive constant $\d>0$ such that, for
any $n$ sufficiently large,
\begin{equation}\label{eq:>d'}
\irn f(u_n) \ge \d.
\end{equation}
\end{enumerate}
\end{lemma}

\begin{proof}
1.\quad For $n$ sufficiently large, by ({\bf f4}), we have
\begin{align*}
b  + \|u_n\| &\ge I_n(u_n) - \frac 1\a \langle I_n'(u_n),u_n
\rangle
\\
&= \left(\frac 12 - \frac 1\a \right)\| u_n\|^2 +\irn V_n(x)\left(
\frac 1\a f'(u_n)u_n - f(u_n) \right)
\\
&\ge \left(\frac 12 - \frac 1\a \right)\| u_n\|^2.
\end{align*}
\\
2. \quad Suppose, by contradiction, that inequality \eqref{eq:br}
does not hold. Then, for any $R>0,$ we should have
\[
\liminf_n \sup_{y\in \RN}  \int_{B_R(y)}|u_n|^2 \, d x=0.
\]
By \cite[Lemma 2]{ABDF}, up to a subsequence,
\[
\lim_n \irn f(u_n)=0
\]
which, by ({\bf f2}) and ({\bf f3}), implies also
\begin{equation*}
\irn f'(u_n)u_n\to 0.
\end{equation*}
Therefore
\begin{align*}
a+o_n(1) &\le I_n(u_n)-\frac 12 \langle I_n'(u_n),u_n\rangle
\\
&=\irn V_n(x)\left( \frac 12 f'(u_n)u_n - f(u_n) \right)=o_n(1),
\end{align*}
which contradicts $a>0.$
\\
By ({\bf f2}) we get \eqref{eq:>d'}.
\end{proof}

\begin{lemma}\label{le:>d}
Let $u_{n,j} \in \D$, $n\ge 1,$ $j\ge 1,$ such that
$\|u_{n,j}\|\ge C>0$ and
\begin{equation}\label{eq:t-c}
\max_{\t\ge 0}I_n(\t u_{n,j}) \le c_n+\d_j,
\end{equation}
with $\d_j\to0^+.$
Then, there exist a sequence $(y_n)_n \subset
\RN$ and two positive numbers $R,\;\mu
>0$ such that
\begin{equation}\label{eq:br-2}
\liminf_n \int_{B_R(y_n)}|u_n|^2 \, d x >\mu,
\end{equation}
where we have set $u_n:=u_{n,n}.$
\\
In particular, there exists a
positive constant $\d>0$ such that, for any $n$ sufficiently
large,
\begin{equation}\label{eq:>d}
\irn f(u_n) \ge \d.
\end{equation}
\end{lemma}

\begin{proof}
Observe that, for any fixed $u\in \D$, $u\neq 0$, there exists
$\bar \t>0$ such that $I_n( \t u)<0$ for any $\t\ge\bar\t$.
\\
As a consequence, the map $g_{n,u}: [0,1]\to\D$ defined by
$$
g_{n,u}(\t)=\t\bar\t u
$$
is in $\G_n$ (which is defined in a natural way) and
\[
\max_{\t \in [0,1]}I_n(g_{n,u}(\t)) = \max_{\t \ge 0}I_n(\t u).
\]
For any $u_{n,j},$ consider the map $g_{n,j}$ defined as before.
By \eqref{eq:t-c} and \cite[Theorem~4.3]{MW}, there exist two
sequences $(w_{n,j})_{ n, j}\subset\D$ and $(\t_{n,j})_{n,j} \subset
[0,1]$ such that
\begin{align}
&\|w_{n,j} - g_{n,j}(\t_{n,j})\| \le \d_j^{1/2}, \label{eq:t-w1}
\\
&|I_n(w_{n,j}) -c_n| < \d_j,  \nonumber 
\\
&\|I_n'(w_{n,j})\| \le \d_j^{1/2}.  \nonumber 
\end{align}
Now we set $w_n:=w_{n,n}$ and analogously we do for $u_{n,n},$ $\t_{n,n}$ and
$g_{n,n}.$ By definition, for $n\ge 1$, there exists $t_n>0$ such
that $g_n(\t_n)=t_n u_n$.
\\
Since $(w_n)_n$ satisfies the hypotheses of Lemma \ref{le:tea}, it
is bounded and there exist a sequence $(y_n)_n \subset \RN$ and
two positive numbers $R,\;\mu
>0$ such that
\begin{equation*}
\liminf_n \int_{B_R(y_n)}|w_n|^2 \, d x >\mu.
\end{equation*}
Moreover, by \eqref{eq:t-w1}, we have
\[
C t_n \le \| t_n u_n\| \le \|t_n u_n -w_n\| + \|w_n\| \le
h_n^{1/2}  + \|w_n\| \le C',
\]
that is $(t_n)_n$ is bounded.
\\
So \eqref{eq:br-2} follows immediately observing that
\begin{align*}
\mu^{1/2} < \|w_n\|_{L^2(B_R(y_n))} &\le \|w_n - t_n
u_n\|_{L^2(B_R(y_n))} +\|t_n u_n\|_{L^2(B_R(y_n))}
\\
&\le C''\big(\|w_n -t_n u_n\| +  \|u_n\|_{L^2(B_R(y_n))}\big)
\\
&\le C''\big(h_n^{1/2} +  \|u_n\|_{L^2(B_R(y_n))}\big).
\end{align*}
By ({\bf f2}), we get \eqref{eq:>d}.
\end{proof}

Let $\widehat V$ be another potential satisfying ({\bf V1-2}) and assume the following notations:
\[
\hat I(u)=\irn |\n u|^2 -\irn \widehat V(x) f(u)\,d x,
\]
$\widehat\Ne$ is its Nehari manifold and $\hat c=c(\widehat V).$

\begin{lemma}\label{le:t-bdd}
Let $(u_n)_n  \subset \D$ such that $\|u_n\|=1$ and
\begin{equation*}
I(\t(u_n)u_n)= \max_{\t\ge 0}I(\t u_n) \to c(V),\qquad \hbox{ as }n
\to \infty.
\end{equation*}
If $\widehat V$ is another potential satisfying ({\bf V1-2})
(eventually $\widehat V=V$),  then the
sequence $(\hat\t(u_n))_n\subset \R_+$ such that for every $n$
$$
\hat I(\hat \t(u_n)u_n)=\max_{\t\ge 0}\hat I( \t u_n),
$$
possesses a
bounded subsequence in $\R$.
\end{lemma}

\begin{proof}
If, up to a subsequence, for all $n \ge 1$, $\hat\t(u_n)\le 1$,
then we are done. Suppose that $\hat\t(u_n)>1$. Then, for all $n
\ge 1$, by ({\bf f4}), we have
\begin{align*}
[\hat\t(u_n)]^2 \irn |\n u_n|^2 &=\irn \widehat V(x)
f'(\t(u_n)u_n) \hat\t(u_n) u_n
\\
&\ge \a \irn \widehat V(x) f(\hat\t(u_n)u_n)
\\
&\ge \a [\hat\t(u_n)]^\a \irn \widehat V(x) f(u_n).
\end{align*}
Since $\a>2$, the conclusion follows from Lemma \ref{le:>d} and ({\bf V2}).
\end{proof}

\begin{lemma}\label{le:VV}
Let $f$ satisfy $(${\bf f1-4}$)$,
$V$ and $\widehat V$ satisfy $(${\bf V1-2}$)$.
\begin{enumerate}
\item If $V\le \widehat V$, then $c\ge \hat c.$
\item If there exists $\delta>0$ such that $V+\delta\le\widehat V,$ then $c>\hat c.$
\end{enumerate}
\end{lemma}

\begin{proof}
1.\quad For all $u\in\D,$ $u\neq 0$, we have
\begin{align*}
\hat c=\inf_{u\neq 0}\sup_{\t\ge 0}\hat I(\t u)
\le\sup_{\t\ge 0}\hat I(\t u)
\le\sup_{\t\ge 0} I(\t u)
\end{align*}
and then the conclusion.
\\
2.\quad By contradiction, suppose $c=\hat c$ and let
$(u_n)_n\subset\Ne$ be such that
\begin{equation}\label{eq:hatc}
I(u_n)\to \hat c.
\end{equation}
Consider the sequence
$(\hat\t(u_n))_n\subset \R_+$ such that, for every $n,$
$$
\hat I(\hat \t(u_n)u_n)=\max_{\t\ge 0}\hat I( \t u_n).
$$
We have
\begin{align*}
I(u_n)&\ge I(\hat\t(u_n) u_n)=\hat
I(\hat\t(u_n)u_n)+\irn (\widehat V(x) -V(x))f(\hat\t(u_n)u_n)\\
&\ge \hat c +\delta\irn f(\hat\t(u_n)u_n),
\end{align*}
so, by \eqref{eq:hatc}, we deduce that
\begin{equation*}
\irn f(\hat\t(u_n)u_n)\to 0.
\end{equation*}
By ({\bf f2-3}) and ({\bf V2}),
$$
\irn \widehat V(x) f'(\hat\t(u_n)u_n)\,\hat\t(u_n)u_n\to 0,
$$
so, since $\hat\t(u_n)u_n\in\widehat\Ne,$ we conclude that
$$
\|\hat\t(u_n)u_n\|\to0.
$$
This fact contradicts 3 of Lemma \ref{le:N}.
\end{proof}

\begin{lemma}\label{le:vn}
Suppose that $f$ satisfies $(${\bf f1-4}$)$ and $V,$ $V_n$ satisfy
$(${\bf V1-2}$)$, for all $n\ge~\!\!1$.
\\
If $V_n\to V$ in $L^\infty(\RN)$ then $c(V_n)\to
c(V).$
\end{lemma}

\begin{proof}
In this proof we repeat the arguments of \cite[Theorem~3.21]{R},
so we skip some details. It is easy to see that we are reduced to
prove the case $V_n=V+h_n,$ with $h_n\to0.$ We first show
$$
c^+:=\lim_{h_n\to 0^+}c(V+h_n)=c(V).
$$
By Lemma \ref{le:VV} certainly  $c^+\le c(V).$ By contradiction
suppose
\begin{equation}\label{eq:cv}
c^+< c(V).
\end{equation}

Let $\delta_j\to 0^+.$ For every $n,j\ge 1,$ by the definition of
$c_n,$ there exists $u_{n,j}\in\D$ such that $\|u_{n,j}\|=1$ and
$$
\max_{\t\ge 0}I_n(\t u_{n,j})\le c_n +\d_j.
$$
Denoting $u_n=u_{n,n},$ since $\D \hookrightarrow L^{2^*}(\RN)$, we have
\begin{align*}
c(V) &\le \max_{\t\ge 0} I(\t u_n) =I(\t(u_n) u_n)
\\
&=I_n(\t(u_n) u_n) +h_n \irn f(\t(u_n) u_n)
\\
&\le \max_{\t\ge 0}I_n (\t u_n) +h_n \irn f(\t(u_n) u_n)
\\
&\le c_n +\d_n +h_n \irn f(\t(u_n) u_n)
\\
&\le c^+ +\d_n +h_n \|\t(u_n)u_n\|_{L^{2^*}}^{2^*}
\\
&\le c^+ +\d_n + C h_n\big(\t(u_n)\big)^{2^*}.
\end{align*}
By Lemma \ref{le:t-bdd} $(\t(u_n))_n$ is bounded, and then we get
a contradiction with \eqref{eq:cv}.
\\
Now we show
$$
c^-:=\lim_{h_n\to 0^-}c(V+h_n)=c(V).
$$
By Lemma \ref{le:VV} certainly  $c^-\ge c(V).$ By contradiction
suppose
\begin{equation*}
c^-> c(V).
\end{equation*}
Let $\d_n \to 0^+$. For every $n\ge 1,$ by the definition of
$c(V),$ there exists a sequence $(u_{n})_n\subset\D$ such that
$\|u_{n}\|=1$ and
$$
\max_{\t\ge 0}I(\t u_{n})\le c(V) +\d_n.
$$
We have
\begin{align*}
c^- &\le c_n \le \max_{\t\ge 0} I_n(\t u_n) =I_n(\t_n(u_n) u_n)
\\
&=I(\t_n(u_n) u_n) -h_n \irn f(\t_n(u_n) u_n)
\\
&\le \max_{\t\ge 0}I (\t u_n) -h_n \irn f(\t_n(u_n) u_n)
\\
&\le c(V) +\d_n -h_n \irn f(\t_n(u_n) u_n)\\
&\le c(V)+\d_n-C h_n \big(\t_n(u_n)\big)^{2^*}.
\end{align*}
Again, the conclusion follows from Lemma \ref{le:t-bdd}.
\end{proof}

In the sequel we will use the following notations
\begin{align*}
V_0      &=\sup_{x\in\RN}V(x);
\\
V_\infty &=\limsup_{|x|\to\infty}V(x).
\end{align*}
By ({\bf V2}), $V_0,V_\infty\in\R_+.$
Moreover we define
\begin{align*}
I_\infty(u) &:= \frac 12 \irn |\n u|^2 - \irn V_\infty f(u),
\\
\Ne_\infty &:= \left\{u\in \D\setminus \{0\}\;\Big{|} \;\irn |\n u|^2=\irn V_\infty f'(u)u \right\},
\\
c_\infty &:= \inf_{u \in \Ne_\infty} I_\infty(u).
\end{align*}

The following theorem is a crucial step in view of the proof of Theorem~\ref{th:main}.

\begin{theorem}\label{th:c}
Suppose that $(${\bf f1-4}$)$ and $(${\bf V1-2}$)$ hold. Let $\widehat V>0$ be
such that
\[
V_\infty\le \widehat V.
\]
Then either $c$ is a critical value of $I$ or $c\ge \hat c$.
\end{theorem}

\begin{proof}
Suppose
\begin{equation}\label{eq:<V}
V_\infty < \widehat V.
\end{equation}
By Lemma \ref{le:ccc}, there exists a sequence $(u_n)_n$ in $\D$,
such that $\|u_n\|=1$ and
\begin{equation}\label{eq:c}
\max_{\t \ge 0} I(\t u_n) \to c, \qquad \hbox{ as }n \to \infty.
\end{equation}
For any $u_n$, we construct the function $g_n\in\Gamma$ as in the
proof of Lemma \ref{le:>d}. Since for any $n\ge 1$
\[
\max_{\t \in [0,1]}I(g_n(\t)) = \max_{\t \ge 0}I(\t u_n),
\]
by \cite[Theorem 4.3]{MW} there exist a sequence $(w_n)_n$ in
$\D$, $h_n>0$, $h_n \to 0$ and $\t_n \in [0,1]$ such that
\begin{align}
&\|w_n - g_n(\t_n)\| \le h_n^{1/2}, \label{eq:w1}
\\
&|I(w_n) -c| < h_n,  \label{eq:w2}
\\
&\|I'(w_n)\| \le h_n^{1/2}.  \nonumber 
\end{align}
Since $(w_n)_n$ is a (PS)-sequence at level $c$, by Lemma
\ref{le:tea} it is bounded and therefore there exists $w\in \D$
such that, up to a subsequence,
\begin{align}
w_n \rightharpoonup w,& \qquad \hbox{ weakly in }\D,  \nonumber
\\
w_n \to w,& \qquad \hbox{ strongly in }L^p_{loc}(\RN).
\label{eq:loc}
\end{align}
It is easy to see that $w$ is a critical point of $I$, then we
need only to check whether $w\neq 0$.
\\
By \eqref{eq:<V}, there exists $\rho>0$, such that for all $|x|\ge
\rho$ we have $V(x)\le \widehat V$. Then, for all $\a>0$, we get
\begin{align*}
\max_{\t \ge 0} I(\t u_n) \ge & \;I(\a u_n)
\\
= & \;\hat I (\a u_n) +\int_{B_\rho}\left(\widehat V - V(x)\right)
f(\a u_n)
\\
&+\int_{\RN \setminus B_\rho}\left(\widehat V - V(x)\right) f(\a u_n)
\\
\ge &\;\hat I (\a u_n) +\int_{B_\rho}\left(\widehat V - V(x)\right)
f(\a u_n).
\end{align*}
Taking $\a=\hat \t(u_n)$, where $\hat \t(u_n)>0$ is such that
\[
\hat I(\hat \t(u_n)u_n) = \max_{\t \ge 0}\hat I(\t u_n),
\]
by Lemma \ref{le:ccc}, referred to $\hat I$, we get
\begin{equation}\label{eq:VV}
\max_{\t \ge 0} I(\t u_n) \ge \hat c +\int_{B_\rho}\left(\widehat V -
V(x)\right) f(\hat \t(u_n) u_n).
\end{equation}
By Lemma \ref{le:t-bdd}, $(\hat \t(u_n))_n$ is a bounded sequence.
\\
Now, according to the definition of $g_n$, for every $n\ge 1$
consider the number $t_n>0$ such that $g_n(\t_n)=t_n u_n$; by
\eqref{eq:w1}
\begin{equation}\label{eq:lowbound}
\|w_n\|_{L^p(B\rho)} \ge \|t_n u_n\|_{L^p(B\rho)} -\|w_n -t_n u_n
\|_{L^p(B\rho)} \ge \|t_n u_n\|_{L^p(B\rho)} - h_n^{1/2}.
\end{equation}
Observe that $(t_n)_n$ is bounded below by a positive constant;
otherwise, since $(u_n)_n$ is a bounded sequence, $I(t_n u_n)\to
0$ along a subsequence, which contradicts \eqref{eq:w1} and \eqref{eq:w2}.
\\
We consider two possibilities:
\begin{itemize}
\item there exists a positive constant $\g$ such that, for any
$n\ge1,$
\begin{equation}\label{eq:g}
\|u_n\|_{L^p(B\rho)}\ge \g;
\end{equation}
\item up to subsequences,
\begin{equation}\label{eq:nong}
\|u_n\|_{L^p(B\rho)} \to 0,\qquad \hbox{ as }n \to \infty.
\end{equation}
\end{itemize}
If \eqref{eq:g} holds, then from \eqref{eq:lowbound} we deduce
that there exists a positive constant $\g'$ such that
\[
\|w_n\|_{L^p(B\rho)} \ge \g'
\]
and this, by \eqref{eq:loc}, ensures that $w\neq 0$.
\\
Moreover $I(w)=c.$ Indeed, since $w\in\Ne,$ certainly $I(w)\ge c.$
On the other hand, by ({\bf f4}), for any $\rho'>0$
\begin{align*}
I(w_n) - \frac 12\langle I'(w_n),w_n \rangle &= \irn V(x) \left(
\frac{1}{2} f'(w_n)w_n - f(w_n) \right)
\\
&\ge \int_{B_{\rho'}} V(x) \left( \frac{1}{2} f'(w_n)w_n - f(w_n)
\right)
\end{align*}
and then, passing to the limit, by \eqref{eq:loc} and the
arbitrariness of $\rho'$, we have
\[
c \ge \irn V(x) \left( \frac{1}{2} f'(w)w - f(w) \right) = I(w).
\]
Hence $c$ is a critical value for $I$.
\\
Suppose, at contrary, that \eqref{eq:nong} holds. Then, by
\eqref{eq:c}, \eqref{eq:VV}, Lemma \ref{le:t-bdd} and the
continuity of the functional
    $$u\in L^p(B_\rho)\mapsto \int_{B_\rho} f(u),$$
we have that $c \ge \hat c.$
\\
Finally, if
\[
V_\infty= \widehat V,
\]
then the conclusion follows from Lemma \ref{le:vn}, using similar
arguments as in~\cite{R}.
\end{proof}

\begin{theorem}\label{th}
Suppose that $f$ satisfies $(${\bf f1-4}$)$ and $V$ satisfies $(${\bf V1-3}$)$.
Then $c$ is a critical value for $I$.
\end{theorem}

\begin{proof}
We apply Theorem \ref{th:c} for $\widehat V=V_\infty.$
\\
By \cite{BL1} (see also \cite{BM}), there exists $w$, a ground
state solution for the problem
\[
\left\{
\begin{array}{l}
-\Delta u = V_\infty f'(u),\quad \hbox{in }\RN,
\\
u>0,
\\
u \in \D,
\end{array}
\right.
\]
namely, $w\in\Ne_\infty$ and $I_\infty(w)=c_\infty$.
\\
Let $\t(w)>0$ be such that $I(\t(w)w)=\max_{\t \ge 0}I(\t w)$. By
({\bf V3}), we have
\begin{align*}
c_\infty & =I_\infty(w) \ge I_\infty(\t(w) w)
\\
&=I(\t(w) w) + \irn(V(x)-V_\infty) f(\t(w)w)
> c,
\end{align*}
and hence, by Theorem \ref{th:c}, we conclude.
\end{proof}

\begin{proofmain}
By the previous theorem, there exists $u\in \D$ such that $I(u)=c$ and $I'(u)=0$. First of all, we prove that $u$ does not change sign. Suppose by contradiction that $u=u^+ + u^-$,
$u^\pm \neq 0$, where $u^+=\max\{0,u\}$ and $u^-=\min\{0,u\}$. It is easy to see that $u^\pm \in \Ne$, so $I(u^\pm)\ge c$: the contradiction arises observing that $I(u)=I(u^+)+I(u^-)$.
\\
Now, since $f$ is even, we can suppose that $u\ge 0$. By the strong maximum principle, we argue that $u>0$  and so it is a solution to problem \eqref{eq}.
\end{proofmain}

If we look for solutions of the problem
\begin{equation}\tag{${\cal P}_\eps$}
\left\{
\begin{array}{l}
-\eps^2 \Delta u - V(x)f'(u)=0,\quad \hbox{in }\RN,
\\
u>0,
\\
u \in \D,
\end{array}
\right.
\end{equation}
for $\eps>0$ sufficiently small, we can weaken the hypotheses on
$V$, replacing ({\bf V3}) by ({\bf V4}).

By the change of variable $x\mapsto \eps x$, the equation \eqref{eq:eps} can be
reduced to the following one
\begin{equation*}
-\Delta u=V(\eps x)f'(u),
\end{equation*}
whose solutions correspond to the critical points of the
functional defined on $\D$
\begin{equation*}
I_\varepsilon(u)=\frac 1 2 \irn|\nabla u|^2\,dx-\irn V(\varepsilon x) f(u)\,dx
\end{equation*}
restricted on the Nehari manifold
\begin{equation*}
\Ne_\eps:=\left\{u\in\D\setminus \{0\} \: \Big| \: \intrn |\nabla u|^2\,dx=\intrn V(\varepsilon x)
f'(u)u\,dx\right\}.
\end{equation*}
We denote by $c_\eps$ the mountain pass level of the
functional $I_\eps$, namely
\[
c_\eps=\inf_{u\in \Ne_\eps}I_\eps(u).
\]
By means of Theorem \ref{th:c}, we will prove that, for small $\eps$, $c_\eps$ is a critical value for $I_\eps$.

We need two preliminary lemmas.
As in Lemma 3.2 of \cite{BM} we can prove the following
\begin{lemma}\label{le:bound}
There exists $C>0$ such that for all $\eps>0$ and, for all $u\in
\Ne_\eps$, we get $\|u\|\ge C.$
\end{lemma}

Now fix $\eta\in \RN$ and let
\begin{align*}
I_\eta(u)&=\frac 1 2 \intrn|\nabla u|^2\,dx-\intrn V(\eta)
f(u)\,dx,
\\
\Ne_\eta&=\left\{u\in\D\setminus \{0\} \:\Big|\: \irn |\n u|^2\,dx
=\irn V(\eta) f'(u)u\,dx\right\},
\end{align*}
and
$c(\eta):=c(V(\eta))$ be the mountain pass level of $I_\eta.$
Consider $\o^\eta$ a ground state solution of the problem
\[
\left\{
\begin{array}{l}
-\Delta u = V(\eta) f'(u),\quad \hbox{in }\RN,
\\
u>0,
\\
u \in \D,
\end{array}
\right.
\]
for any $\eps>0$ define
$$
\o^\eta_\eps=\o^\eta(\cdot-\eta/\eps)
$$
and let $\t^\eta_\eps>0$ be such that
$\t^\eta_\eps \o^\eta_\eps\in\Ne_\eps.$ The following result holds
\begin{lemma}\label{le:c-eta}
For any $\eta \in \RN$, we get
\begin{equation*}
\lim_{\eps\to 0} I_\eps(\t^\eta_\eps \o^\eta_\eps)=c(\eta).
\end{equation*}
\end{lemma}

\begin{proof}
First we show that $(\t^\eta_\eps)_{\eps>0}$ is bounded. If
$\t^\eta_\eps \le 1$, we are done; otherwise by some computations we
have
\begin{align*}
(\t^\eta_\eps)^2\irn |\n \o^\eta_\eps|^2
&=\irn V(\eps x) f'(\t^\eta_\eps\o_\eps^\eta)\t^\eta_\eps\o_\eps^\eta
\ge \alpha\irn V(\eps x)f(\t^\eta_\eps\o^\eta_\eps)
\\
&\ge C\a(\t^\eta_\eps)^\a\irn f(\o^\eta_\eps)
\end{align*}
and then, by a change of variable,
\begin{equation*}
(\t^\eta_\eps)^2\irn |\n \o^\eta|^2 \ge C\a(\t^\eta_\eps)^\a\irn f(\o^\eta),
\end{equation*}
from which we deduce that $(\t^\eta_\eps)_{\eps>0}$ is bounded.
\\
Let $\t^\eta\ge 0$ be such that, up to a subsequence, $\t_\eps^\eta\to\t^\eta$, as $\eps\to 0$. Since $\t_\eps^\eta \o^\eta_\eps\in\Ne_\eps$, by Lemma \ref{le:bound} we have that
\begin{equation*}
\t^\eta_\eps\|\o^\eta\|=\t_\eps^\eta\|\o^\eta_\eps\|=\|\t_\eps^\eta\o^\eta_\eps\|\ge C.
\end{equation*}
and then $\t^\eta\neq 0.$ We prove that $\t^\eta=1.$
\\
Indeed
\begin{align*}
(\t^\eta_\eps)^2\irn V(\eta)f'(\o^\eta)\o^\eta
& = (\t_\eps^\eta)^2\irn|\n\o^\eta|^2
= (\t_\eps^\eta)^2\irn|\n\o_\eps^\eta|^2
\\
& = \irn V(\eps x)f'(\t_\eps^\eta\o^\eta_\eps)\t_\eps^\eta\o_\eps^\eta
\\
& = \irn V(\eps x + \eta) f'(\t_\eps^\eta\o^\eta)\t_\eps^\eta\o^\eta.
\end{align*}
Letting $\eps \to 0$ and using the Lebesgue's theorem,
\begin{equation*}
(\t^\eta)^2\irn V(\eta) f'(\o^\eta) \o^\eta=\irn V(\eta) f'(\t^\eta\o^\eta)\t^\eta\o^\eta,
\end{equation*}
so
\begin{equation*}
\irn \big( f'(\t^\eta\o^\eta)\o^\eta-f'(\o^\eta)\t^\eta\o^\eta \big)=0.
\end{equation*}
Since for any $z\in\R$, $z\neq 0$, the function
$$t>0 \mapsto \frac{f'(t z)z}{t}-f'(z) z$$
vanishes only for $t=1$, we deduce that $\t^\eta=1.$
\\
In conclusion
\begin{align*}
I_\eps(\t^\eta_\eps \o_\eps^\eta)
& = \frac{(\t_\eps^\eta)^2}{2}\irn|\n \o^\eta_\eps|^2
-\irn V(\eps x) f(\t^\eta_\eps \o_\eps^\eta)
\\
& = \frac{(\t_\eps^\eta)^2}{2}\irn|\n \o^\eta|^2
-\irn V(\eps x+\eta) f(\t^\eta_\eps \o^\eta)
\to I_\eta(\o^\eta)=c(\eta),
\end{align*}
and the proof is complete.
\end{proof}

Arguing as in the proof of Theorem \ref{th:main}, Theorem \ref{th:maineps} is an immediate consequence of the following
\begin{theorem}\label{th:2}
Suppose that $f$ satisfies $(${\bf f1-4}$)$ and $V$ satisfies $(${\bf V1-2}$)$ and $(${\bf V4}$)$. Then there exists $\bar \eps >0$ such that
for any $\eps\in (0,\bar\eps)$, $c_\eps$ is a critical value for
$I_\eps$.
\end{theorem}

\begin{proof}
Suppose by contradiction that for any $\bar\eps>0$ there exists
$\eps<\bar\eps$ such that $c_\eps$ is not a critical value for
$I_\eps$. Then, by Theorem \ref{th:c}, there exists a sequence
$\eps_n\searrow 0^+$ such that $(c_{\eps_n})_n$ is bounded from
below by $c_\infty.$
\\
By ({\bf V4}) there exists $\eta \in \RN$ such that
$V(\eta)>V_\infty$, so, by 2 of Lemma \ref{le:VV},
\[
c(\eta)<c_\infty\le c_{\eps_n}.
\]
On the other side, by Lemma \ref{le:c-eta},
we know that
\[
c_{\eps_n}\le I_{\eps_n}(\t^\eta_{\eps_n} \o_{\eps_n}^\eta) \to
c(\eta)
\]
and so, for $\eps_n$ sufficiently small, we get a contradiction.
\end{proof}

\section{A multiplicity result}\label{sec:mult}

This section is devoted to the proof of Theorem \ref{th:main2}. In
view of this, from now on we assume that all the hypotheses of
Theorem \ref{th:main2} hold.

Set
\begin{equation*}
c_0:=\inf_{\eta\in\R^N} c_\eta.
\end{equation*}
By Lemmas \ref{le:VV} and \ref{le:vn}, we have that
$$
c_0=c(V_0)=\inf_{u\in\Ne_0} I_0(u),
$$
where
\begin{align*}
I_0(u)&:=\frac 1 2 \irn|\n u|^2-\irn V_0 f(u),
\\
\Ne_0&:=\left\{u\in\D\setminus\{0\}\;\Big|\;\irn|\n u|^2=\irn V_0 f'(u)u\right\}.
\end{align*}
As a consequence,
$$
M:=\Big\{\eta\in \RN \;\big|\; V(\eta)=\max_{\xi\in \RN}V(\xi)\Big\}
= \Big\{\eta\in\R^N \;\big|\;  c_\eta=c_0\Big\};
$$
moreover, Lemma \ref{le:VV} and ({\bf V4}) imply that $M$ is
compact and
\begin{equation}\label{eq:0-inf}
c_0<c_\infty.
\end{equation}
For all $a\in\R$ and $\eps>0$, we define $I_\eps^a:=\{u\in \Ne_\eps\mid
I_\eps(u)\le a\}$.

To prove Theorem \ref{th:main2} we will refer to the following abstract multiplicity
theorem (see \cite{MW})

\begin{theorem}\label{th:abs}
Let $\M$ be a $C^{1,1}$ complete Riemannian manifold
modeled on an Hilbert space and $J$ be a $C^1$ functional on $\M$
bounded from below. If there exists $b>\inf_{\M}J$ such that
$J$ satisfies the Palais-Smale
condition on the sublevel $J^{-1}(-\infty,b)$, then for any
noncritical level $a,$ with $a<b,$ there exist at least
$\cat_{J^a}(J^a)$ critical points of $J$ in $J^a,$ where
$J^a:=\{u\in\M\mid J(u)\le a\}.$
\end{theorem}
So, in order to solve \eqref{eq:eps}, we need to study the
topology of the sublevels of the functional ${I_\eps}|_{\Ne_\eps}$, which is positive by ({\bf f4}). In
particular, we will compare the topology of the sublevels of $I_\eps$
with that of $M$ using the following lemma, which is a consequence
of the definitions of category and homotopic equivalence (we refer
to \cite{BC} for more details)

\begin{lemma}\label{le:cat}
Let $\eps>0$, $a\in\R$ and $\g>0.$ If there exist
$\psi:M\to I_\eps^a$ and $\beta:I_\eps^a\to M_\g$ two continuous maps
such that
$\beta\circ\psi$ is homotopically equivalent to the
embedding $j:M\to M_\g,$ then
$\cat_{I_\eps^a}(I_\eps^a)\ge\cat_{M_\g}(M).$
\end{lemma}

Taking these two results into account, the proof of
Theorem \ref{th:main2} can be divided in two steps: the study of
the topology and the study of the compactness of the sublevels.
\\
The subsection \ref{sec:top} will be devoted to the construction
of the maps $\psi$ and $\beta$ in such a way we can relate the
topology of a suitable sublevel of $I_\eps|_{\Ne_\eps}$ with that
of $M.$
\\
In subsection \ref{sec:comp} we prove the compactness of the
Palais-Smale sequences in a suitable sublevel of
$I_\eps|_{\Ne_\eps}$, which is guaranteed by assumption $({\bf
V4})$ and a concentration-compactness argument.
\\
Finally, in subsection \ref{sec:proof} we give the proof of Theorem \ref{th:main2}.

\subsection{The topology of the sublevels}\label{sec:top}

Fix $\g>0$. For any $\eps>0$ define the map $\b_\eps : \D \setminus \{0\} \to \RN$ as
\begin{equation*}
\beta_\varepsilon(u)
=\frac{\irn|\n u|^2\chi(\varepsilon x)\,dx}{\irn|\n u|^2\,dx},
\quad \hbox{for all }u\in\D\setminus \{0\},
\end{equation*}
where $\chi:\RN \to \RN$ is defined as
\begin{equation*}
\chi(x):=\left\{
\begin{array}{ll}
x & \hbox{ if }|x|\le\rho,
\\
\rho\frac{x}{|x|} & \hbox{ if }|x|>\rho,
\end{array}
\right.
\end{equation*}
with $\rho>0$ such that $M_\g\subset B_\rho=\{x\in\RN\mid
|x|<\rho\}.$

It is easy to see that for any $\eps>0$ the map $\beta_\eps$ is
continuous.
\begin{lemma}\label{le:beta}
For any $u\in\D\setminus \{0\},$ $\eps>0,$ $\eta\in M,$ denote by
$$
u_{\eps,\eta}:x\in\RN\mapsto u(x-\eta/\eps)\in\R.
$$
Then
\begin{equation}\label{eq:bar}
\lim_{\eps\to 0}\beta_\eps(u_{\eps,\eta})=\eta,
\end{equation}
uniformly in $M.$
\end{lemma}

\begin{proof}
By some computations
\begin{align*}
\beta_\eps(u_{\eps,\eta})
& = \frac{\irn|\n u_{\eps,\eta}|^2\chi(\varepsilon x)\,dx}
{\irn|\n u_{\eps,\eta}|^2\,dx}
= \frac{\irn|\nabla u|^2\chi(\varepsilon x+\eta)\,dx}
{\irn|\n u|^2\,dx}
\\
& = \eta + \frac{\irn|\nabla u|^2\big(\chi(\varepsilon x+\eta)-\eta\big)\,dx}
{\irn|\n u|^2\,dx}.
\end{align*}
So
\begin{equation*}
\lim_{\eps\to 0}\sup_{\eta\in M}\big|\beta_\eps(u_{\eps,\eta})-\eta\big|\le
\lim_{\eps\to 0}\sup_{\eta\in M}
\frac{\irn|\nabla u|^2\big|\chi(\varepsilon x+\eta)-\eta\big|\,dx}
{\irn|\n u|^2\,dx}=0
\end{equation*}
since, by the compactness of $M,$ for any $\d>0$ there
exist $r,\bar \eps>0$ such that for all $\eta\in M$ and for all $\eps\in (0,\bar\eps)$
$$
\irn|\nabla u|^2\big|\chi(\varepsilon x+\eta)-\eta\big|\,dx
\le \eps r\int_{B_r}|\nabla u|^2 \,dx+2\rho\int_{B_r^c}|\n u|^2\,dx\le \d.
$$
\end{proof}

Now we introduce a technical lemma which describes a sort of compactness for any sequence $(u_n)_n$ such that for all $n \ge 1$, $u_n \in  \Ne_{\eps_n}$ and $I_{\eps_n}(u_n) \to c_0$. Observe that such sequences exist by the definition of $c_0$ and by Lemma \ref{le:c-eta}. In the proof, we will follow an idea of \cite{AF}.

\begin{lemma}\label{le:claim}
Let $\eps_n \to 0^+$, as $n\to \infty$, and, for all $n\ge 1$,
$u_n\in \Ne_{\eps_n}$ such that
\begin{equation}\label{eq:limc0}
\lim_n I_{\eps_n}(u_n)=c_0.
\end{equation}
Then there exists a sequence $(\eta_n)_n$ in $\RN$, $\eta \in M$
and $v\in \D$, such that
\begin{enumerate}
\item $\eta_n \to \eta$, as $n\to \infty$;
\item $v_n:=u_n \left(\cdot +\eta_n/\eps_n\right) \to v$ in $\D$, as $n\to
\infty$.
\end{enumerate}
\end{lemma}

\begin{proof}
Since for any $n\ge 1 $ $u_n\in \Ne_{\eps_n},$ by ({\bf f4}) we
have
\begin{align*}
c_0 + o_n(1)
&=\left(\frac 12 - \frac 1\a \right)\| u_n\|^2
+\irn V(\eps x)\left( \frac 1\a f'(u_n)u_n - f(u_n) \right)
\\
&\ge \left(\frac 12 - \frac 1\a \right)\| u_n\|^2
\end{align*}
and then $(u_n)_n$ is bounded in $\D.$ Using \cite[Lemma 2]{ABDF},
by similar arguments as in 2 of Lemma \ref{le:tea}, we can prove
that there exists a sequence $(\xi_n)_n\subset\RN$ and two
positive constants $R,$ $\mu>0$ such that for any $n$ large enough
\begin{equation}\label{eq:concentration}
\int_{B_R(\xi_n)}|u_n|^2\ge \mu.
\end{equation}
Define $v_n:=u_n(\cdot+\xi_n),$ $\eta_n:=\eps_n\xi_n$ and $\t_n>0$
such that, for any $n\ge 1,$ $\t_n v_n\in\Ne_0.$
\\
\\
{\sc Claim 1:} there exists a positive constant $C$ such that $
(\t_n)_n\subset [C, 1].$
\\
\\
Since $(v_n)_n$ is bounded, by 3 of Lemma \ref{le:N} certainly
$(\t_n)_n$ is bounded below by some $C>0.$ Moreover, since for any
$n\ge 1$
\begin{align*}
\t_n^2\irn|\n v_n|^2&=\irn V_0f'(\t_n v_n)\t_n v_n
\\
\noalign{\hbox{and}}
\irn |\n v_n|^2 &= \irn V(\eps_n x+\eta_n)f'(v_n)v_n,
\end{align*}
we have
\begin{equation*}
\irn V_0 f'(\t_n v_n)\t_n v_n
=\t_n^2\irn V(\eps_n x+\eta_n)f'(v_n)v_n
\le \t_n^2\irn V_0f'(v_n)v_n,
\end{equation*}
that is
\begin{equation*}
\irn
V_0\left(\frac{f'(\t_n v_n)v_n}{\t_n}-f'(v_n)v_n\right)\le 0.
\end{equation*}
We conclude the proof of the claim just observing that for any
$z\in\R$, $z\neq 0$, the function
    \begin{equation}\label{eq:functionn}
        t>0 \mapsto\frac{f'(t z)z}{t}-f'(z) z
    \end{equation}
is non positive if and only if $t\le 1.$
\\
\\
{\sc Claim 2:} $I_0(\t_n v_n)\to c_0$.
\\
\\
Since $(\t_n v_n)_n\subset\Ne_0,$ we have
\begin{align*}
c_0
&\le I_0(\t_nv_n)
=\frac 1 2 \irn|\n(\t_n v_n)|^2
-\irn V_0 f(\t_n v_n)
\\
&\le\frac 1 2 \irn|\n(\t_n v_n)|^2
-\irn V(\eps_n x + \eta_n) f(\t_n v_n)
\\
&=\irn V(\eps_n x +\eta_n)\left(\frac{\t_n^2}{2} f'(v_n)v_n-f(\t_n v_n)\right)
\\
&\le\irn V(\eps_n x+\eta_n)\left(\frac{1}{2} f'(v_n)v_n-f(v_n)\right)
=I_{\eps_n}(u_n)\to c_0,
\end{align*}
where we have used the fact that, for any $z\in\R,$ $z\neq 0,$ the
function
\begin{equation}\label{eq:cresce}
t\in[0,1]\mapsto (t^2/2)f'(z)z-f(t z)
\end{equation}
is increasing. Now define $w_n=\t_n v_n.$
\\
\\
{\sc Claim 3:} $(w_n)_n$ converges strongly in $\D$ to some $w$
which is a ground state solution of the problem
\begin{equation}\label{eq:P0}
\left\{
\begin{array}{ll}
-\Delta u - V_0f'(u)=0, &\hbox{in }\RN,
\\
u>0,
\\
u \in \D.
\end{array}
\right.
\end{equation}
\\
By Claim 2 and taking \cite[Theorem 8.5]{Wi} into account, we can
suppose that the sequence $(w_n)_n$ satisfies the (PS)-condition
for the functional ${I_0}|_{\Ne_0}$; by this assumption, it can be
proved (see e.g. \cite{CL2}) that $(w_n)_n$ is also a
(PS)-sequence for the unconstrained functional. By Claim 1, the
sequence $(w_n)_n$ is bounded and then there exists $w\in\D$ such
that, up to a subsequence,
\begin{eqnarray}
w_n&\rightharpoonup& w \;\hbox{weakly in }\D,\label{eq:wd1}
\\
w_n&\to& w\;\hbox{in }L^s(B), \hbox{ with $B\subset \RN$, bounded,
and }1\le s<2^*.\label{eq:slp1}
\end{eqnarray}
Observe that $w\in\Ne_0;$ indeed, by \eqref{eq:concentration},
Claim 1 and \eqref{eq:slp1} we deduce that $w\neq 0$, while from
\eqref{eq:wd1} and \eqref{eq:slp1} it follows that
$I'_0(w)=0.$
\\
So, for any $\d>0$, there exists $r'=r'(\d)>0$ such that
\begin{align*}
c_0&\le I_0(w) = \irn V_0\left(\frac 1 2 f'(w)w-f(w)\right)
\\
&\le \int_{B_{r'}}V_0\left(\frac 1 2 f'(w)w - f(w)\right) + \d
\\
&= \lim_n \int_{B_{r'}}V_0 \left(\frac 1 2 f'(w_n)w_n-f(w_n)\right)+\d,
\end{align*}
from which we deduce that
\begin{align}
\limsup_n \int_{B^c_{r'}}  V_0 \Big(\frac \alpha 2 -1\Big)f(w_n)
&\le \lim_n\int_{B^c_{r'}}V_0\left(\frac 1 2 f'(w_n)w_n-f(w_n)\right)\nonumber
\\
&= c_0- \lim_n \int_{B_{r'}}V_0\left(\frac 1 2 f'(w_n)w_n-f(w_n)\right) \nonumber
\\
&\le\d.     \label{eq:small1}
\end{align}
By ({\bf f2}) and ({\bf f3}) and \eqref{eq:small1} we deduce that,
for any $\d>0$, there exists $r''=r''(\d)>0$ such that
\begin{equation*}
\limsup_n\int_{B^c_{r''}}V_0 f'(w_n)w_n\le \d,
\end{equation*}
therefore for any $\d>0$
\begin{align}
\irn|\n w|^2
&\le \liminf_n\irn|\n w_n|^2
\le \limsup_n\irn V_0 f'(w_n)w_n\nonumber
\\
&=\lim_n\int_{B_{r''}} V_0 f'(w_n)w_n
+\limsup_n\int_{B^c_{r''}} V_0 f'(w_n)w_n\nonumber
\\
&\le \int_{B_{r''}}V_0 f'(w)w + \d
\le\irn|\n w|^2 +\d.\label{eq:chain}
\end{align}
By \eqref{eq:wd1} and \eqref{eq:chain} it follows that, up to a
subsequence, $w_n\to w$ in $\D$ and then $w$ is a ground state
solution of \eqref{eq:P0}.
\\
\\
{\sc Claim 4:} $(\eta_n)_n$ converges to some $\eta\in M$ and
$(v_n)_n$ converges strongly in $\D$ to a ground state solution of
\eqref{eq:P0}.
\\
\\
First observe that, by Claim 1, up to a subsequence $(\t_n)_n$
converges to some $\t_0>0.$ Therefore, by Claim 3, there exists a
subsequence (identically relabeled) of $(v_n)_n$ and
$v\in\D\setminus\{0\}$ such that $v_n\to v$ in $\D.$
\\
There are two possibilities:
\begin{enumerate}
\item $|\eta_n|\to +\infty$;
\item up to a subsequence, there exists $\eta \in \RN$ such that $\eta_n\to\eta\in\RN.$
\end{enumerate}
Suppose that $|\eta_n|\to +\infty$ as $n \to \infty$.
\\
For any fixed $\d>0$, let $r=r(\d)>0$ be such that
\[
\int_{B_r^c} f(v_n) < \d.
\]
Since $V_\infty=\limsup_{|x|\to \infty}V(x)$, for $n$ sufficiently
large, and for all $x\in B_r$, we get
\[
V(\eps_n x +\eta_n)\le V_\infty +\d.
\]
Therefore, for $n$ large
\begin{multline*}
\irn V(\eps_n x +\eta_n)f(v_n)  =\int_{B_r} V(\eps_n x
+\eta_n)f(v_n) +O(\d)
\\
\le \int_{B_r} (V_\infty +\d)f(v_n) +O(\d)\le \irn (V_\infty
+\d)f(v_n) +O(\d).
\end{multline*}
Passing to the limit and by the arbitrariness of $\d>0$
\begin{equation}\label{eq:limsup}
\limsup_n \irn V(\eps_n x +\eta_n)f(v_n) \le \irn V_\infty f(v).
\end{equation}
Let $\t_\infty>0$ be such that $\t_\infty v \in \Ne_{\infty}$,
namely
\[
\t_\infty^2\irn |\n v|^2 =\irn V_\infty f'(\t_\infty v)\t_\infty
v.
\]
By \eqref{eq:limsup} and since $u_n \in \Ne_{\eps_n}$, we have
\begin{align*}
c_\infty &\le \frac{\t_\infty^2}{2}\irn |\n v|^2 - \irn V_\infty
f(\t_\infty v)
\\
&\le \frac{\t_\infty^2}{2} \lim_n \irn |\n v_n|^2 - \limsup_n \irn
V(\eps_n x +\eta_n) f(\t_\infty v_n)
\\
&\le \liminf_n \left(\frac{\t_\infty^2}{2} \irn |\n v_n|^2 - \irn
V(\eps_n x +\eta_n) f(\t_\infty v_n)\right)
\\
&= \liminf_n I_{\eps_n} (\t_\infty u_n) \le \liminf_n I_{\eps_n}
(u_n)= c_0,
\end{align*}
and we get a contradiction with \eqref{eq:0-inf}.
\\
So, up to subsequences, there exists $\eta\in \RN$ such that
$\eta_n \to \eta$, as $n \to \infty$. By Lebesgue theorem,
\begin{align*}
\irn V(\eps_n x +\eta_n) f(v_n)
&\to \irn V(\eta) f(v)
\\
\noalign{\hbox{and}}
\irn V(\eps_n x +\eta_n) f'(v_n)v_n
&\to \irn V(\eta) f'(v)v,
\end{align*}
from which we deduce
\begin{equation*}
I_\eta(v)=c_0\;{\hbox{ and }}\;v\in\Ne_\eta,
\end{equation*}
that is $c_0=c_\eta$ and $\eta\in M.$
\end{proof}

\begin{theorem}\label{th:inclusion}
Let $\d_n \to 0^+$, as $n \to \infty$. Then, for every $\g>0$, there exists $(\bar\eps_n)_n$, $\bar\eps_n\to 0^+$, such
that, for $n$ sufficiently large and for every $\eps\in (0,\bar \eps_n)$, $I_\eps^{c_0+\d_n}\neq \emptyset$ and  $\beta_\eps(I_\eps^{c_0+\d_n})\subset M_\g.$
\end{theorem}

\begin{proof}
By Lemma \ref{le:c-eta}, certainly for any $n\ge 1$
$I_\eps^{c_0+\d_n}\neq\emptyset$ for small $\eps.$ Now suppose by
contradiction that there exists $\g>0$ and $\eps_n \to 0^+$ such
that for any $n\ge 1$ there exists $u_n\in I_{\eps_n}^{c_0+\d_n}$
and
\begin{equation}\label{eq:dist}
\dist(\b_{\eps_n}(u_n),M)>\g.
\end{equation}
Since by Lemma \ref{le:VV} $c_0\le c_{\eps_n}$ for any  $n\ge 1,$
we have
    $$\lim_n I_{\eps_n}(u_n)=c_0.$$
By Lemma \ref{le:claim}, there exists a sequence $(\eta_n)_n$ in
$\RN$, $\eta \in M$ and $v\in \D$, such that $\eta_n \to \eta$ and
$v_n:=u_n \left(\cdot +\eta_n/\eps_n\right) \to v$, in $\D$, as
$n\to \infty$. This implies that $(\eta_n)_n \subset M_\g$.
\\
We claim that
\[
\lim_n \irn |\n v_n|^2 \chi (\eps_n x  +\eta_n) = \irn |\n v|^2
\eta.
\]
In fact, for any
$\d>0$, there exists $r=r(\d)$, such that, for $n$ sufficiently
large,
\[
\int_{B_r^c} |\n v_n|^2 \le \d,\qquad
\int_{B_r^c} |\n v|^2 \le \d,
\]
hence,
\begin{align*}
&\left| \irn |\n v_n|^2 \chi (\eps_n x +\eta_n)
- \irn |\n v|^2 \eta \right|
\\
& \qquad\quad\le \left| \int_{B_r^c} |\n v_n|^2 \chi (\eps_n x +\eta_n)\right|
+\left| \int_{B_r^c} |\n v|^2 \eta \right|
\\
&\qquad\qquad+\left| \int_{B_r} (|\n v_n|^2 -|\n v|^2) \chi (\eps_n x +\eta_n)\right|
\\
&\qquad\qquad+\left| \int_{B_r} |\n v|^2 (\chi (\eps_n x+\eta_n)
-\eta)\right|=C \d.
\end{align*}
This implies that
\[
\b_{\eps_n}(u_n) =\frac{\irn |\n u_n|^2 \chi (\eps_n x)}{\irn |\n
u_n|^2} =\frac{\irn |\n v_n|^2 \chi (\eps_n x +\eta_n)}{\irn |\n
v_n|^2 } \to \eta
\]
which contradicts \eqref{eq:dist}.
\end{proof}

Let $\o$ be a ground state solution of the problem
\begin{equation*}
\left\{
\begin{array}{ll}
-\Delta u=V_0f'(u), & \hbox{in }\RN,
\\
u>0,
\\
u\in \D.
\end{array}
\right.
\end{equation*}
For any
$\eta\in M$ and $\eps>0$ define the new function
$$
\o^\eta_\eps=\o(\cdot-\eta/\eps)
$$
and let $\t^\eta_\eps>0$ be such that
$\t^\eta_\eps\o^\eta_\eps\in\Ne_\eps.$ We set
\begin{equation*}
\Phi_\eps:\eta\in M\mapsto \t^\eta_\eps\o^\eta_\eps\in\Ne_\eps.
\end{equation*}
By \cite {BM}, $\Phi_\eps$ is continuous. Moreover, arguing as in Lemma \ref{le:c-eta}, we can prove the following result
\begin{theorem}\label{th:conc0}
Uniformly for $\eta\in M$
\begin{equation*}
\lim_{\eps\to 0} I_\eps(\Phi_\eps(\eta))=c_0.
\end{equation*}
\end{theorem}

Combining the results of Lemma \ref{le:beta}, Theorems
\ref{th:inclusion} and \ref{th:conc0} we get the following
\begin{theorem}\label{th:cat}
Let $\d_n\to 0^+$, as $n\to \infty$. Then, for every $\g>0$, there exists $(\bar\eps_n)_n$, $\bar\eps_n\to 0^+$, such
that for $n$ sufficiently large and for every $\eps\in (0,\bar \eps_n)$ we have
$$
\cat_{\tilde I_\eps^n} \tilde I_\eps^n \ge
\cat_{M_\g}M,
$$
where $\tilde I_\eps^n:=I^{c_0+\d_n}_\eps.$
\end{theorem}

\begin{proof}
Let $\d_n\to 0^+$, as $n\to \infty$, and $\g>0$. According to Theorem \ref{th:conc0}, there exists $(\bar\eps_n')_n$ such that for every $\eps\in (0,\bar \eps_n')$
\begin{equation}\label{eq:Phi2}
\Phi_\eps:\eta\in M\mapsto \Phi_\eps(\eta)\in I_\eps^{c_0+\d_n}.
\end{equation}
By Theorem \ref{th:inclusion}, there exists $(\bar\eps_n'')_n$, $\bar\eps_n''\to 0^+$, such
that, for $n$ sufficiently large and for every $\eps\in (0,\bar \eps_n'')$:
\begin{equation}\label{eq:beta2}
\beta_\eps:u\in I_\eps^{c_0+\d_n}\mapsto
\beta_\eps (u)\in M_\g.
\end{equation}
These last two formulas hold simultaneously for any $\eps\in (0,\bar \eps_n)$, where $\bar\eps_n=\min\{\bar\eps_n',\bar\eps_n''\}$.
\\
Moreover using Lemma \ref{le:beta} we have that, uniformly
for $\eta\in M$
\begin{equation*}
\lim_{\eps\to 0}\beta_\eps(\Phi_\eps(\eta))=\eta.
\end{equation*}
So for every $\eps>0$ sufficiently small, the map $\beta_{\eps}\circ\Phi_{\eps}$ is
homotopically equivalent to the canonical injection
$j:M\to M_\g$. By \eqref{eq:Phi2}, \eqref{eq:beta2} and Lemma
\ref{le:cat} we get the conclusion.
\end{proof}

\subsection{The compactness on the sublevels}\label{sec:comp}

This section is completely devoted to the study of the compactness
properties of the Palais Smale sequences. In particular, in view
of Theorem \ref{th:abs} and of the topological considerations in
the previous section, we are interested in investigating the
compactness properties of the sublevels of the type $I^a_\eps$
with $a>c_0.$ The following result has been obtained by similar
arguments as in \cite{AF}.

\begin{lemma}\label{le:comp}
Let $(v_n)_n\subset\D$ and $d < c_\infty$  be such that
\begin{equation}\label{eq:quasiPS}
I_\eps(v_n)\to d, \quad\quad \langle I_\eps'(v_n), v_n \rangle \to 0.
\end{equation}
If
\begin{equation*}
v_n\rightharpoonup 0\quad\hbox{in }\D,
\end{equation*}
then $v_n\to 0$ in $\D.$
\end{lemma}
\begin{proof}
    Let $(v_n)_n$ satisfy \eqref{eq:quasiPS} with
    $d<c_\infty,$ and assume that $v_n\rightharpoonup0.$
    \\
    We show that $d=0.$ Indeed, since $\langle
    I'_\eps(v_n),v_n\rangle\to0,$ we have
        \begin{equation}\label{eq:PScond}
            \irn|\n v_n|^2=\irn V(\eps x)f'(v_n)v_n +o_n(1),
        \end{equation}
    and then, by ({\bf f4}),
        \begin{align*}
            d & = I_\eps(v_n) -\frac 1 2 \langle I'_\eps(v_n),v_n\rangle + o_n(1)\\
            &= \irn V(\eps x) \left(\frac 1 2 f'(v_n)v_n -
            f(v_n)\right) + o_n(1)\ge o_n(1),
        \end{align*}
    from which we deduce that $d\ge 0.$ Now suppose by contradiction
    that $d>0.$
    It is easy to see that the conclusions of Lemma \ref{le:tea} hold also for $(v_n)_n$, so there exist a sequence $(y_n)_n \subset \RN$
    and three positive numbers $R,\;\mu,\;\d>0$ such that
        \begin{equation}\label{eq:brr}
            \liminf_n \int_{B_R(y_n)}|v_n|^2 \, d x >\mu,
        \end{equation}
    and
        \begin{equation}\label{eq:>dd'}
            \irn f(v_n) \ge \d.
        \end{equation}
    Let $(\t_n)_n\subset\R_+$ be such that
    $\t_n v_n\in\Ne_\infty.$ We prove
    that $(\t_n)_n$ is bounded. If $\t_n\le 1$ we are done; otherwise,
    since by ({\bf f4})
    $$
        \t_n^2\irn|\n v_n|^2 =\irn V_{\infty}f'(\t_nv_n)\t_n v_n \ge
        \a \t_n^{\a}\irn V_{\infty} f(v_n),
    $$
    the conclusion follows by the boundedness of $(v_n)_n$ in $\D$,
    ({\bf f2-3}) and \eqref{eq:>dd'}.
    \\
    We are going to prove by contradiction that $\liminf_n\t_n
    \le 1.$ Define
    $\tilde v_n:=v_n(\cdot+y_n)$ and let $\rho>0$ and $R'>0$ such
    that
        \begin{equation*}
            V(\eps x)\le V_\infty+\rho,\quad \forall\, |x|\ge R'.
        \end{equation*}
    We have that, for any $(w_n)_n\subset\D$ such that
    $w_n\rightharpoonup 0$ weakly in $\D$
        \begin{align}\label{eq:chainchain1}
            \irn V(\eps
            x) f'(w_n)w_n \nonumber & = \int_{B_{R'}}\!\! V(\eps
            x) f'(w_n)w_n+\int_{B_{R'}^c}\!\! V(\eps
            x) f'(w_n)w_n\nonumber\\
            &\le o_n(1)+\int_{B_{R'}^c}\!\! (V_\infty+\rho)f'(w_n)w_n\nonumber\\
            &\le o_n(1) + O(\rho) + \irn V_\infty f'(w_n)w_n,
        \end{align}
    and, analogously,
        \begin{equation}\label{eq:chainchain2}
            \irn V(\eps
            x) f(w_n)\le o_n(1) + O(\rho) + \irn V_\infty
            f(w_n).
        \end{equation}
    Since $\t_n v_n\in \Ne_\infty,$ by \eqref{eq:PScond} and \eqref{eq:chainchain1}
    we have
        \begin{align}\label{eq:chainchain}
            \irn V_\infty f'(\t_nv_n)\t_nv_n +o_n(1) & = \t_n^2\irn V(\eps
            x) f'(v_n)v_n \nonumber\\
            &\le o_n(1) + O(\rho) + \t_n^2\irn V_\infty f'(v_n)v_n.
        \end{align}
    If we suppose that $\liminf_n\t_n > 1,$ then by \eqref{eq:functionn}
    and \eqref{eq:chainchain}
        \begin{align*}
            \int_{B_R}\left(\frac{f'(\t_n \tilde v_n)
            \tilde v_n}{\t_n}-f'(\tilde v_n)\tilde v_n\right)&\le
            \irn \left(\frac{f'(\t_n \tilde v_n)
            \tilde v_n}{\t_n}-f'(\tilde v_n)\tilde v_n\right)\\
            &\le o_n(1)+O(\rho).
        \end{align*}
On the other hand, by the boundedness of $(\tilde v_n)_n$ and of
$(\t_n)_n$, from \eqref{eq:brr} we deduce that, up to a
subsequence,
\begin{equation*}
\int_{B_R}\left(\frac{f'(\t_n\tilde v_n)\tilde v_n}{\t_n} -
f'(\tilde v_n)\tilde v_n\right) \to C >0.
\end{equation*}
Then, up to a subsequence, one of the following two possibilities
holds:
\begin{description}
\item{$i)$}  $\forall n\ge 1: \t_n\le 1$,\\
\item{$ii)$} $\forall n\ge 1: \t_n\ge 1$ and $\lim_n\t_n=1.$
\end{description}
If $i)$ holds, then by \eqref{eq:cresce}, \eqref{eq:PScond} and
\eqref{eq:chainchain2} we have
\begin{align*}
c_\infty &\le \frac 1 2 \irn |\n (\t_n v_n)|^2 -\irn V_\infty
f(\t_n v_n)
\\
&\le \irn V(\eps x)\left(\frac {\t_n^2} 2 f'( v_n) v_n - f(\t_n
v_n)\right) + o_n(1) + O(\rho)
\\
&\le \irn V(\eps x)\left(\frac 1 2 f'( v_n) v_n - f(v_n)\right) +
o_n(1) + O(\rho)
\\
&=I_\eps (v_n) + o_n(1) + O(\rho);
\end{align*}
if $ii)$ holds, then, by \eqref{eq:chainchain2},
\begin{align*}
I_\eps(v_n)- c_\infty &\ge \frac 1 2 \irn|\n v_n|^2 - \irn V(\eps
x) f(v_n)
\\
&\quad - \frac {\t_n^2} 2 \irn|\n v_n|^2 + \irn V_\infty f(\t_n
v_n)
\\
& \ge \frac{1-\t_n^2}{2}\irn |\n v_n|^2 +\irn V(\eps x)\big(f(\t_n
v_n)-f(v_n)\big)
\\
&\quad + o_n(1) + O(\rho)\\
 &\ge o_n(1) + O(\rho).
\end{align*}
Both in the first and in the second case we can conclude that
\begin{align*}
c_\infty \le I_\eps(v_n)+O(\rho) + o_n(1) = d+O(\rho) + o_n(1),
\end{align*}
and then, letting $n$ go to $\infty$ and taking $\rho$ smaller and
smaller, we deduce $c_\infty\le d$ which contradicts our
hypothesis.
\\
So we have proved $d=0,$ that is
    \begin{equation*}
        \irn V(\eps x) \left(\frac 1 2 f'(v_n)v_n -
        f(v_n)\right) \to 0.
    \end{equation*}
By ({\bf f4}) we deduce that
    \begin{equation*}
        \irn V(\eps x) f(v_n)\to0
    \end{equation*}
and then, by ({\bf f2}), ({\bf f3}) and \eqref{eq:PScond},
    \begin{equation*}
        \irn |\n v_n|^2 = \irn V(\eps x) f' (v_n) v_n + o_n(1)\to
        0,
    \end{equation*}
and we are done.
\end{proof}
\begin{theorem}\label{th:comp}
For any $\eps>0$ small enough, the sublevel $I_\eps^{c_\infty}$ is
nonempty and, moreover, $I_\eps|_{\Ne_\eps}$ satisfies the
(PS)-condition in the strip $[c_\eps,c_\infty).$
\end{theorem}

\begin{proof}
First observe that, by Theorem \ref{th:conc0} and hypothesis ({\bf
V4}), for $\eps$ small enough the sublevel $I_\eps^{c_\infty}$ is
nonempty.
\\
Now, let $(u_n)_n\subset \Ne_\eps$ be a Palais Smale sequence at
the level $\lambda< c_\infty,$ namely
\begin{align}
&I_\eps(u_n)    =  \lambda + o_n(1),          \label{eq:PS1}
\\
&{I_\eps'}|_{\Ne_\eps}(u_n)  = o_n(1).          \label{eq:PS2}
\end{align}
Actually $(u_n)_n$ is a (PS)-sequence for the unconstrained
functional, namely
\begin{equation}\label{eq:PS3}
\lim_n\sup_{\substack{ v\in D^{1,2} \\ \|v\|=1}} \langle
I'_\eps(u_n),v\rangle =0.
\end{equation}
By Lemma \ref{le:tea}, the sequence $(u_n)_n$ is bounded in $\D,$
and therefore there exists $u\in\D$ such that, up to a
subsequence,
\begin{eqnarray}
u_n&\rightharpoonup& u \;\hbox{weakly in }\D,\label{eq:wd}
\\
u_n&\to& u\;\hbox{in }L^s(B), \hbox{ with $B\subset \RN$, bounded, and }1\le s<2^*.\label{eq:slp}
\end{eqnarray}
By \eqref{eq:PS2}, \eqref{eq:wd} and \eqref{eq:slp}
$I'_\eps(u)=0,$ so
    \begin{equation}\label{eq:Iu}
        I_\eps(u) = I_\eps(u) - \frac 1 2 \langle I'_\eps(u),u\rangle =
        \irn  V(\eps x) \bigg[\frac 1 2 f'(u)u -
        f(u)\bigg]
        \ge 0.
    \end{equation}
We set $v_n=u_n-u,$ so that our aim is to prove that $v_n\to 0.$
We show that $(v_n)_n$ satisfies all the hypotheses of Lemma
\ref{le:comp}.
\\
By \cite[Lemma 2.8]{BM} and \eqref{eq:Iu}
    \begin{equation*}
        I_\eps(v_n)=I_\eps(u_n)-I_\eps(u) +
        o_n(1)=\l-I_\eps(u)+o_n(1)\to \l-I_\eps(u) < c_\infty.
    \end{equation*}
Moreover, again by \cite[Lemma 2.8]{BM}, we infer that
\[
\langle I'_\eps(v_n),v_n\rangle
=\langle I'_\eps(u_n),u_n\rangle
-\langle I'_\eps(u),u\rangle
+o_n(1)
=o_n(1),
\]
and so we are done.
\end{proof}

\subsection{Proof of Theorem \ref{th:main2}}\label{sec:proof}

Let $\g>0$ and fix $\d_n\to 0^+$.
\\
Since $c_0<c_\infty$, for $n$ sufficiently large, $I_\eps^{c_0+\d_n}\subset I_\eps^{c_\infty}$ and then Theorem \ref{th:comp} implies that (PS)-condition holds in $I_\eps^{c_0 +\d_n}$, for small $\eps$. Therefore, applying Theorem \ref{th:abs} to our case, there exists at least $\cat_{I_\eps^{c_0 +\d_n}}(I_\eps^{c_0 +\d_n})$ critical points of the functional $I_\eps$. Now by Theorem \ref{th:cat}, up to take smaller $\eps$ and greater $n$, we find at least $\cat_{M_\g}(M)$ critical points of $I_\eps$ with energy less or equal to $c_0+\d_n$. We need only to prove that such solutions are strictly positive. First we show that they do not change sign. Otherwise, we would have $u\in \D$ a critical point of $I_\eps$,
\begin{equation}\label{eq:pos}
I_\eps(u)\le  c_0+\d_n,
\end{equation}
such that $u=u^+ + u^-$, $u^\pm \neq 0$, where $u^+=\max\{0,u\}$ and $u^-=\min\{0,u\}$. Since $u^\pm \in \Ne_\eps$, then $I_\eps(u^\pm)\ge c_\eps \ge c_0$. But $I_\eps(u)=I_\eps(u^+)+I_\eps(u^-)\ge 2 c_0$ which contradicts \eqref{eq:pos}.
\\
Now, since $f$ is even, we can suppose that all these solutions are nonnegative. Actually, by the strong maximum principle, we argue that they are positive.

\end{document}